\documentclass[11pt]{article}

\makeatletter
\@addtoreset{equation}{section}
\makeatother

\newtheorem{remark}{Remark}
\newtheorem{theorem}{Theorem}

\input epsf

\usepackage{amssymb}
\usepackage{amsmath}
\usepackage{epsfig}

\def\tfrac#1#2{{\frac{\lower.6ex
\hbox{$\scriptstyle#1$}}
{\raise.7ex
\hbox{$\scriptstyle#2$}}}}

\def\sign{{\rm sign}}

\def\bigO{{\cal O}}

\def\calC{{\cal C}}

\def\calL{{\cal L}}
\def\phase{{\rm ph}}
\def\dsp#1{\displaystyle#1}

\def\binomial#1#2{
\renewcommand{\arraystretch}{1.0}
\left(
\begin{array}{c} 
\hskip-5pt#1\\
\hskip-5pt#2
\end{array}
\hskip-5pt\right)}

\def\F11#1#2#3{
{}_1F_1\left(
\begin{array}{c}
\begin{array}{c}\hskip-10pt#1\end{array}\\
\begin{array}{c}\hskip-10pt #2\end{array}
\end{array}
\hskip-8pt;\,#3
\right)}

\def\protectbold#1{\protect{\boldmath{$#1$}}}

\begin{document}

\title{Large Degree Asymptotics of Generalized Bessel Polynomials}

\author{Jos\'e Luis L\'opez \\
    Departamento de Ingenieria Mat\'ematica e Inform\'atica,\\
    Universidad P\'ublica de Navarra, 31006-Pamplona, Spain\\
       \and
    Nico M. Temme\\
    Centrum Wiskunde \& Informatica, \\
    Science Park 123, 1098 XG Amsterdam, The Netherlands. \\
     { \small e-mail: {\tt
    jl.lopez@unavarra.es, 
    nicot@cwi.nl}}
    }

\date{\today}
\maketitle

\begin{abstract}
\noindent
Asymptotic expansions are given for large values of $n$ of the generalized Bessel polynomials $Y_n^\mu(z)$. The analysis is based on integrals that follow from the generating functions of the polynomials. A new simple expansion is given that is valid outside a compact neighborhood of the origin in the $z-$plane.  New forms of expansions in terms of elementary functions valid  in sectors not containing the turning points $z=\pm i/n$ are derived, and a new expansion in terms of modified Bessel functions is given. Earlier asymptotic expansions of the generalized Bessel polynomials by Wong and Zhang (1997) and Dunster (2001) are discussed.

\end{abstract}

\vskip 0.8cm \noindent
{\small
\par\noindent
Keywords \& Phrases:
generalized Bessel polynomials, 
asymptotic expansions, 
Bessel functions, 
Kummer functions.
}
\section{Introduction}\label{sec:intro}

Generalized  Bessel polynomials of degree
$n$, complex order $\mu$ and complex argument $z$, denoted by $Y_n^\mu(z)$, have been introduced in \cite{Krall:1949:ANC}, and can be defined by their generating function. We have \cite {Grosswald:1978:BPS}:
\begin{equation}\label{Intro1}
\left(\frac2{1+\sqrt{{1-2zw}}}\right)^\mu\,\frac{e^{2w/(1+\sqrt{{1-2zw}})}}{\sqrt{{1-2zw}}}=
\sum_{n=0}^\infty \frac{Y_n^{\mu}(z)}{ n!}w^n,
\quad \vert 2zw\vert<1,
\end{equation}
with special values $Y_n^{\mu}(0)=1$, $n=0,1,2,\ldots$.

The generalized Bessel polynomials are important in certain
problems of mathematical physics. For a historical survey and discussion of many interesting
properties, we refer to  \cite {Grosswald:1978:BPS}.

In \cite{Dunster:2001:UAE}  and \cite{Wong:1997:AEG} detailed contributions on asymptotic approximations are given on the generalized Bessel polynomials. Wong and Zhang use integral representations and Dunster's approach is based on a differential equation. Our approach also uses integrals as starting point. 

Our approach is different from that of Wong and Zhang \cite{Wong:1997:AEG}. We give expansions that are similar to those of the modified Bessel functions, and the expansions reduce to these expansions when $\mu=0$. For the expansions  in terms of elementary functions we give a simple description for the domains of validity.

In \S\ref{sec:Besalt} we start with a simple expansion valid outside a fixed neighborhood of the origin.  In \S\ref{sec:elem} we give expansions that hold uniformly with respect to $z$ inside the sectors $\vert\phase\,\pm z\vert\le \frac12\pi-\delta$, where $\delta$ is a small positive number. At the end of this section we compare our results with those in \cite{Dunster:2001:UAE}  and \cite{Wong:1997:AEG}. For complementary sectors (and extensions of these) we give an expansion in terms of the modified Bessel function $K_{n+\frac12}(z)$ and its derivative. In fact, this is an expansion in terms  of $Y_n^{0}(z)$, and this expansion holds for all $z$. In the Appendix  \S\ref{sec:modBes} we summarize the expansions of the modified Bessel functions. These expansions play a role when comparing the results for the generalized Bessel polynomials with those for $\mu=0$.


\section{Representations and relations with Bessel and Kummer functions}\label{sec:Rels}
For $\mu=0$ the generalized Bessel polynomials become well-known polynomials that occur in representations of Bessel functions of fractional order. We have in terms of the modified $K-$Bessel function 
\cite{Olver:2010:BFS}\footnote{http://dlmf.nist.gov/10.49.E12}:
\begin{equation}\label{ybes}
Y_n^{0}(z)=\sqrt{\frac{2}{\pi z}}e^{1/z}K_{n+\frac12}(1/z)=
\sum_{k=0}^n \binomial{n}{k}\,(n+1)_k\,\left(\tfrac12 z\right)^k,
\end{equation}
where  $(p)_k$ is the Pochhammer symbol defined by
\begin{equation}\label{i4}
(p)_0=1,\quad (p)_k=\frac{\Gamma(p+k)}{\Gamma(p)},\quad k\ge1.
\end{equation}
For $Y_n^{\mu}(z)$ an explicit formula reads \cite{Krall:1949:ANC}
\begin{equation}\label{ydef}
Y_n^{\mu}(z)=\sum_{k=0}^n\binomial{n}{ k}\,(n+\mu+1)_k\,\left(\tfrac12z\right)^k.
\end{equation}

The simple integral representation
\begin{equation}\label{yint}
Y_n^{\mu}(z)=\frac1{\Gamma(n+\mu+1)}\int_0^\infty t^{n+\mu}\left(1+\tfrac12zt\right)^n e^{-t}\,dt
\end{equation}
gives the representation in \eqref{ydef} by expanding $(1+\frac12zt)^n$ in powers of $z$. 

From the generating function in \eqref{Intro1} we have the Cauchy integral representation
 \begin{equation}\label{Besint}
Y_n^{\mu}(z)=\frac{n!}{2\pi i}\int_{\calC}\left(\frac2{1+\sqrt{{1-2zw}}}\right)^\mu\,\frac{e^{2w/(1+\sqrt{{1-2zw}})}}{\sqrt{{1-2zw}}}\,\frac{dw}{w^{n+1}},
\end{equation}
where $\calC$ is a circle with radius smaller than $1/\vert 2z\vert$. 

When $z=0$ all polynomials reduce to unity, the first few polynomials being
\renewcommand{\arraystretch}{1.75}
\begin{equation}\label{y012}
\begin{array}{l}
\dsp{Y_0^{\mu}(z)=1,\quad Y_1^{\mu}(z)=1+\tfrac12(\mu+2)z,}\\ 
\dsp{Y_2^{\mu}(z)=1+(\mu+3)z+\tfrac14(\mu+3)(\mu+4)z^2.}
\end{array}
\end{equation}
\renewcommand{\arraystretch}{1}%
More values can be obtained from the recurrence relation
\begin{equation}\label{yrecn}
A_nY_{n+2}^{\mu}(z)=B_nY_{n+1}^{\mu}(z)+C_nY_{n}^{\mu}(z),
\end{equation}
where
\renewcommand{\arraystretch}{1.75}
\begin{equation}\label{yrecc}
\begin{array}{l}
\dsp{A_n=2(2n+\mu+2)(n+\mu+2),}\\ 
\dsp{B_n=(2n+\mu+3)(2\mu+z(2n+\mu+4)(2n+\mu+2)),}\\
\dsp{C_n=2(n+1)(2n+\mu+4).}
\end{array}
\end{equation}
\renewcommand{\arraystretch}{1.0}%
There is also a recursion with respect to $\mu$:
\begin{equation}\label{yrecmu}
(n+\mu+2)Y_{n}^{\mu+2}(z)=(2n+\mu+2-2/z)Y_{n}^{\mu+1}(z)+(2/z)Y_{n}^{\mu}(z),
\end{equation}
and for the derivative we have
\begin{equation}\label{ydiff}
\frac{d}{dz}Y_{n}^{\mu}(z)=\tfrac12n(n+\mu+1)Y_{n-1}^{\mu+2}(z)=
\frac{n+\mu+1}{z}\left(Y_{n}^{\mu+1}(z)-Y_{n}^{\mu}(z)\right).
\end{equation}
These relations are  special cases of  known analytic continuation formula for the Kummer functions. 
See, for example, \cite{Olde:2010:CHF} and \cite[p.~13]{Slater:1960:CHF}. They follow also from \eqref{yint} by integrating by parts.

The relations with the Kummer functions  are
\renewcommand{\arraystretch}{2.0}
\begin{equation}\label{ykum}
\begin{array}{@{}r@{\;}c@{\;}l@{}}
Y_n^{\mu}(1/z)&=&\dsp{(2z)^{n+\mu+1}U(n+\mu+1,2n+\mu+2,2z),}\\
&=&\dsp{(2z)^{-n}\frac{\Gamma(2n+\mu+1)}{\Gamma(n+\mu+1)}{}_1F_1(-n;-2n-\mu;2z).}
\end{array}
\end{equation}
For $\Re z<0$ it is convenient to have the representation
\begin{equation}\label{ykumneg}
Y_n^{\mu}(-1/z)=F_n^{\mu}(1/z)+U_n^{\mu}(1/z),
\end{equation}
where
\begin{equation}\label{ykumnegFU}
\begin{array}{l}
F_n^{\mu}(1/z)=\dsp{\frac{n!\,(2z)^{n+\mu+1}e^{-2z}}{\Gamma(2n+\mu+2)}{}_1F_1(n+1;2n+\mu+2;2z),}\\
U_n^{\mu}(1/z)=\dsp{\frac{(-1)^n n!\,(2z)^{n+\mu+1}e^{-2z}}{\Gamma(n+\mu+1)}U(n+1,2n+\mu+2,2z).}
\end{array}
\end{equation}
\renewcommand{\arraystretch}{1.0}%

For $\mu=0$ we have
\begin{equation}\label{zn06}
F_n^{0}(z)=\sqrt{2\pi z}e^{-z} I_{n+\frac12}(1/z), \quad 
U_n^{0}(z)=(-1)^n\sqrt{\frac{2z}{\pi }}e^{-z}K_{n+\frac12}(z),
\end{equation}
and this corresponds to the relation
\begin{equation}\label{zn05}
Y_n^{0}(-z)=\sqrt{\frac{2}{\pi z}}e^{-1/z}\left((-1)^nK_{n+\frac12}(1/z)+\pi I_{n+\frac12}(1/z)\right).
\end{equation}

\section{A simple expansion}\label{sec:Besalt}
The interesting region  in the $z$ plane for uniform asymptotic expansions is a neighborhood of size $\bigO(1/n)$ of the origin, where the zeros appear. For $z$ outside a fixed neighborhood a simple  expansion will be derived. 

First we mention
\begin{equation}\label{Ba01}
\begin{array}{l}
\dsp{Y_n^{\mu}(z)=2^{\mu+\frac12}\left(\frac{2nz}{e}\right)^n\,e^{1/z} \ \times}\\
\quad\quad\quad\quad
\dsp{\left(1-\frac{1+6\mu(\mu+1+2z^{-1})+6z^{-2}}{24n}+\bigO\left(1/n^{2}\right)\right).}
\end{array}
\end{equation}
This is derived in  \cite{Docev:1962:OGB} and  mentioned in \cite[p.~124]{Grosswald:1978:BPS} and \cite{Wong:1997:AEG}. More terms in this expansion can be obtained, for example by using the Cauchy integral given in \eqref{Besint}.

In this section we derive a simple asymptotic expansion related to the result in \eqref{Ba01} by expanding part of the integrand in \eqref{Besint} in powers of $W=\sqrt{1-2zw}$. First we notice that the main asymptotic contributions from the contour integral in \eqref{Besint} come from the singular point $w=1/(2z)$, and when  $w\sim1/(2z)$ the quantity $W$ is small.

We have
\begin{equation}\label{Ba02}
\left(\frac2{1+\sqrt{{1-2zw}}}\right)^\mu\,e^{2w/(1+\sqrt{1-2zw})}=2^\mu e^{1/z}(1+W)^{-\mu}e^{-W/z}
\end{equation}
and we expand  for $m=0,1,2,\ldots\,$
\begin{equation}\label{Ba03}
(1+W)^{-\mu}e^{-W/z}=\sum_{k=0}^{m-1} L_k^{-\mu-k}(1/z)W^k +W^m U_m(W).
\end{equation}
The appearance of the Laguerre polynomials becomes clear when expanding both the exponential and binomial and by comparing the coefficients with the representation 
\begin{equation}\label{Ba04}
L_n^{\alpha}(x)=\sum_{m=0}^n \binomial{n+\alpha}{n-m}\frac{(-x)^m}{m!}.
\end{equation}
Introducing this expansion in \eqref{Besint} we find
\begin{equation}\label{Ba05}
Y_n^{\mu}(z)= n!\,2^\mu e^{1/z} \sum_{k=0}^{m-1} L_k^{-\mu-k}(1/z)\Phi_k^{(n)}+n!\,R_m(n),
\end{equation}
where
\begin{equation}\label{Ba06}
\Phi_k^{(n)}=\frac{1}{2\pi i}\int_{\calC}(1-2zw)^{(k-1)/2}\frac{dw}{w^{n+1}}=(2z)^n\frac{\left(\tfrac12-\tfrac12k\right)_n}{n!},
\end{equation}
and
\begin{equation}\label{Ba07}
R_m(n)=\frac{2^\mu e^{1/z}}{ 2\pi i}\int_{\calC}(1-2wz)^{(m-1)/2}U_m(W)\,\frac{dw}{ w^{n+1}},
\end{equation}
with $\calC$ a circle with radius less than $1/\vert 2z\vert$. 

After the change of variable $w=(1-t/n)/(2z)$ we have
\begin{equation}\label{Ba08}
R_m(n)=\frac{(2z)^n2^\mu e^{1/z}}{ 2\pi i\,n^{(m+1)/2}}\int_{\calC}t^{(m-1)/2}U_m\left(\sqrt{t/n}\right)\,\frac{dt}{ (1-t/n)^{n+1}}.
\end{equation}
The function $U_m(W)$ is  analytic  in $\vert W\vert<1$ and ${\cal O}(1)$ as $W\to 0$. This means that $\vert U_m(W)\vert<C_m$ with $C_m$ a positive constant, on and inside the path $\calC$ of integration in \eqref{Besint}. Hence, 
$\vert U_m(\sqrt{t/n})\vert<C_m$ on $\calC$ (indeed,  the path $\calC$ has been modified after the change of variable, but we can set it equal to the previous path). Also, $(1-t/n)^{-(n+1)}$ is bounded on $\calC$ (and converges to $e^t$ for large  $n$). Therefore, the above integral is $\bigO(1)$ as $n\to \infty$. Hence, 
\begin{equation}\label{Ba09}
R_m(n)={\bigO}\left(\frac{(2z)^n}{ n^{(m+1)/2}}\right),\quad n\to \infty, \  z\ne0,
\end{equation}
which is comparable with the large $n$ behavior of  $\Phi_m^{(n)}$. 

From the integral in \eqref{Ba06} it easily follows that $\Phi_{2k+1}^{(n)}=0$ for  $k=0, 1, 2, \ldots, n$, and we see that, when $n$ is large, only the even terms in the series give contributions. We notice that $\Phi_{2k}^{(n)}$ constitute an asymptotic sequence. This follows from 
\begin{equation}\label{Ba10}
\frac{\Phi_{2k+2}^{(n)}}{\Phi_{2k}^{(n)}}=\frac{\Gamma(-\frac12-k+n)}{\Gamma(\frac12-k+n)}=\frac{1}{-\frac12-k+n}=\bigO\left(n^{-1}\right),\quad n\to\infty.
\end{equation}

We can collect the results of the section as follows.
\begin{theorem}\label{thm1}
For $n\to \infty$ we have the asymptotic expansion
\begin{equation}\label{Ba11}
Y_n^{\mu}(z)\sim(2z)^n\,2^\mu e^{1/z} \sum_{k=0}^\infty L_k^{-\mu-k}(1/z)\left(\tfrac12-\tfrac12k\right)_n,
\end{equation}
which holds uniformly for $\vert z\vert \ge z_0$, where $z_0$ is a positive constant.
\end{theorem}

\begin{table}
\caption{Relative errors $\delta$ in the asymptotic expansion  in \eqref{Ba11} with terms up and including $k=20$ for $\mu=17/4$  and several values of $z$  and $n$. 
\label{Atab}}
\begin{center}
\begin{tabular}{@{}lrcccrcccc@{}}
\hline
   & &  $z=10^j$ & & &  &$z=-10^j$ & \\
   \hline
 &
$j\ $ & 
$Y_n^{\mu}(z)$&
{$\delta$} &
\quad &
{$j\ $} &
$Y_n^{\mu}(z)$&
{$\delta$} 
   \\ \hline
$ n=50$
& -1 &  0.4232e034 &  0.16e-03 &  & -1 &  0.1961e026 &  0.26e-11 \\ 
&  0 &  0.1211e081 &  0.38e-07 &  &  0 &  0.1778e080 &  0.62e-08 \\ 
&  1 &  0.5131e130 &  0.17e-07 &  &  1 &  0.4235e130 &  0.14e-07 \\ 
&  2 &  0.4707e180 &  0.16e-07 &  &  2 &  0.4617e180 &  0.15e-07 \\ 
&  3 &  0.4666e230 &  0.15e-07 &  &  3 &  0.4657e230 &  0.15e-07 \\ 
  \hline
$n=100$
& -1 &  0.1681e093 &  0.30e-07 &  &-1 &  0.5251e084 &  0.68e-15 \\ 
&  0 &  0.3190e189 &  0.10e-10 &  &  0 &  0.4501e188 &  0.18e-11 \\ 
&  1 &  0.1325e289 &  0.47e-11 &  &  1 &  0.1089e289 &  0.39e-11 \\ 
&  2 &  0.1213e389 &  0.43e-11 &  &  2 &  0.1189e389 &  0.42e-11 \\ 
&  3 &  0.1202e489 &  0.43e-11 &  &  3 &  0.1200e489 &  0.42e-11 \\ 
\hline
\end{tabular}
\end{center}
\end{table}

In Table~\ref{Atab} we give the relative errors $\delta$ when we use the expansion in \eqref{Ba11} with terms up and including $k=20$, for $\mu=17/4$  and several values of $z$  and $n$. We see a quite uniform error with respect to $z$, except when $z=\pm\frac{1}{10}$.

\begin{remark}\label{Barem1}
{\rm
If we wish we can expand the Pochhammer symbols in \eqref{Ba11} for large $n$ and rearrange the series. In that way we can obtain an expansion of $Y_n^{\mu}(z)$ in negative powers of $n$, and this expansion is comparable with an expansion of which the first terms are given in 
\eqref{Ba01}.
}
\end{remark}

\begin{remark}\label{Barem2}
{\rm
In  \eqref{Ba11} we expand the generalized Bessel polynomials $Y_n^\mu(z)$ in terms of another set of polynomials, the generalized Laguerre polynomials. Because the degree of these polynomials does not depend on the large parameter, they can be evaluated much easier than the polynomials $Y_n^\mu(z)$. In fact, to compute the Laguerre polynomials we can use a recurrence relation, which follows from differentiating \eqref{Ba03} with respect to $W$. Let $c_k=L_k^{-\mu-k}(1/z)$, then
\begin{equation}\label{Ba12}
z(k+1)c_{k+1}=-(\mu z+kz+1)c_k-c_{k-1}, \quad k=1, 2, 3,\ldots,
\end{equation}
with initial values $c_0=1$, $c_1=-(\mu z+1)/z$.
}
\end{remark}

\begin{remark}\label{Barem3}
{\rm
For general values of $\mu$ the expansion in \eqref{Ba11} is not convergent, but for $\mu=0,-1,-2,\ldots$ it is. For example, for $\mu=0$ a relation for the $K-$Bessel function should arise. We have from \eqref{Ba03} $L_k^{-k}(1/z)=(-1)^k/(k!\,z^k)$, which gives the convergent expansion
\begin{equation}\label{Ba13}
Y_n^{0}(z)=(2z)^n e^{1/z} \sum_{k=0}^\infty\frac{(-1)^k}{k!\,z^k}\left(\tfrac12-\tfrac12k\right)_n.
\end{equation}
Summing the series, separating the terms with even and odd $k$, we obtain
\begin{equation}\label{Ba14}
Y_n^{0}(z)=(-1)^ne^{1/z} \sqrt{\frac{\pi}{2z}}\left(I_{-n-\frac12}(1/z)-I_{n+\frac12}(1/z)\right),\end{equation}
and by using a well-known relation between the modified Bessel functions the representation in \eqref{ybes} arises.
}
\end{remark}

\begin{remark}\label{Barem4}
{\rm
The expansion in  \eqref{Ba11} is simpler than those of the following sections: it is easier to obtain and the coefficients are easily computed.
When $z$ is small the expansion in \eqref{Ba11} breaks down. The Laguerre polynomials ($k\ge1$) are not bounded, although the factor $(2z)^n$ in front of the expansion has some control. But the main concern is the exponential factor $e^{1/z}$, which has an essential singularity at $z=0$. Recall that the polynomials $Y_n^{\mu}(z)$ all tend to unity when $z\to0$. 
As mentioned in Theorem~\ref{thm1}, for the expansion in  \eqref{Ba11}  we have to exclude a fixed neighborhood of the point $z=0$.
}
\end{remark}

\section{Expansions  in terms of elementary functions}\label{sec:elem}
By using saddle point methods we obtain expansions that hold uniformly inside sectors $\vert\phase\,\pm z\vert\le\frac12\pi-\delta$, where $\delta$ is a fixed small positive number.

For $Y_n^{\mu}(z)$ we take $\nu=n+\frac12$ as the large parameter. This gives a suitable identification of the results with those for the Bessel function $K_{n+\frac12}(z)$ when $\mu=0$; see also \eqref{ybes}. In addition we replace the argument $z$ of $Y_n^{\mu}(z)$ by $1/(\nu z)$ (observe that in \cite{Dunster:2001:UAE} $\nu$ is also the large parameter, and the Bessel polynomial is considered with reversed argument). 

Because for both cases $\vert\phase\,\pm z\vert\le\frac12\pi-\delta$ the derivation of the asymptotic expansion is very similar we first summarize the results in the following two subsections, and in \S\ref{sec:deriv} we give the details of the analysis.

\subsection{Expansion holding for \protectbold{\vert\phase\,z\vert<\frac12\pi}}\label{sec:Beszpos}

\begin{theorem}\label{thm2}
For large values of $n$ we have the expansion
\begin{equation}\label{Ynp15}
Y_n^{\mu}(\zeta)\sim\frac{\left(1-z+\sqrt{1+z^2}\right)^{\mu}\sqrt{z}}{(1+z^2)^{\frac14}}e^{\nu z-\nu \eta}
\sum_{k=0}^\infty \frac{A_{k}(\mu,z)} {\nu^k},
\end{equation}
and the expansion holds uniformly inside the sector $\vert\phase\,z\vert\le\frac12\pi-\delta$. Here, $\delta$ is a small positive constant, $\nu=n+\frac12$, $\zeta=1/(\nu z)$, $A_0(\mu,z)=1$, 
 \begin{equation}\label{Ynp23}
A_1(\mu,z)=\frac{t(5t^2-3)}{24}-\frac{\mu t^2(z+1)}{4}+\frac{\mu^2(tz-1)}{4},
\end{equation}
and
\begin{equation}\label{etat0}
t=\frac{1}{\sqrt{1+z^2}}, \qquad
\eta=\sqrt{1+z^2}+\log\frac{z}{1+\sqrt{1+z^2}}.
\end{equation}

\end{theorem}

For $\mu=0$ the coefficients $A_{k}(\mu,z)$ reduce to those in the expansion in \eqref{Knuzuniform}, that is, $A_{k}(0,z)=(-1)^k u_k(t)$, and
\begin{equation}\label{Ynp17}
Y_n^{0}(\zeta)=\sqrt{\frac{2\nu z}{\pi}}e^{\nu z}K_{\nu}(\nu z)\sim\frac{\sqrt{z}}{(1+z^2)^{\frac14}}e^{\nu z-\nu \eta}
\sum_{k=0}^\infty \frac{(-1)^k u_k(t)} {\nu^k},
\end{equation}
which indeed gives the expansion in  \eqref{Knuzuniform}

For $\zeta=0$ all Bessel polynomials $Y_n^{\mu}(\zeta)$ reduce to unity. We have as  $\zeta\to0$:
 \begin{equation}\label{Ynp24}
z\to\infty, \  t\to0, \  zt\to1, \  z-\eta\to0,  .
\end{equation}
As a consequence, $A_1(\mu,z)\to0$ as $\zeta\to0$. In fact all coefficients $A_k(\mu,z)$ with $k\ge1$ vanish as $\zeta\to0$, and both sides of \eqref{Ynp15} reduce to unity.  

Recall that the simple expansion in \S\ref{sec:Besalt} is no longer valid when the argument of the Bessel polynomials approaches the origin.

\subsection{Expansions holding for \protectbold{\vert\phase(-z)\vert<\frac12\pi}}\label{sec:Beszneg}
In this case we write (see \eqref{ykumneg} and \eqref{ykumnegFU})
\begin{equation}\label{zn02n}
Y_n^{\mu}(-1/z)=F_n^{\mu}(1/z)+U_n^{\mu}(1/z).
\end{equation}

We have the following results.
\begin{theorem}\label{thm3}
For large values of $n$ we have the expansions
\begin{equation}\label{zn13}
U_n^{\mu}(\zeta)\sim(-1)^n\frac{\left(1+z+\sqrt{1+z^2}\right)^{\mu}\sqrt{z}}{(1+z^2)^{\frac14}}e^{-\nu z-\nu \eta}
\sum_{k=0}^\infty \frac{B_{k}(\mu,z)} {\nu^k},
\end{equation}
 \begin{equation}\label{zn27}
F_n^{\mu}(\zeta)\sim
\frac{\left(1+z-\sqrt{1+z^2}\right)^\mu\sqrt{z}}{(1+z^2)^{1/4}}e^{-\nu z +\nu\eta}\sum_{k=0}^\infty\frac{C_k(\mu,z)}{\nu^k},
\end{equation}
and the expansions hold uniformly inside the sector $\vert\phase\,z\vert\le\frac12\pi-\delta$. Here, 
$B_{0}(\mu,z)=1$, $C_0(\mu,z)=1$, 
 \begin{equation}\label{zn15}
B_1(\mu,z)=\frac{t(5t^2-3)}{24}+\frac{\mu t^2(z-1)}{4}-\frac{\mu^2(zt+1)}{4},
\end{equation}
and
 \begin{equation}\label{zn29}
C_1(\mu,\zeta)=-\frac{t(5t^2-3)}{24}+\frac{\mu t^2(z-1)}{4}+\frac{\mu^2(zt-1)}{4}.
\end{equation}
The quantities $t$,  $\zeta$,  $\nu$ and $\eta$ are as in Theorem~\ref{thm2}. 
\end{theorem}

For $\mu=0$ the expansions reduce to those for the modified Bessel functions mentioned in \eqref{zn06}.

\subsection{Integral representations}\label{sec:intrep}
For deriving the asymptotic expansions we introduce the integrals
\renewcommand{\arraystretch}{1.75}
\begin{equation}\label{intrep1}
\begin{array}{l}
\dsp{P_\nu^\mu(z)=\int_0^\infty p_\mu(s) e^{-\nu\phi(s)}\,ds,}\\
\dsp{Q_\nu^\mu(z)=\frac{1}{2\pi i}\int_\calL  q_\mu(s) e^{\nu\phi(s)}\,ds,}
\end{array}
\renewcommand{\arraystretch}{1}
\end{equation}
where $\nu>0$ and
\begin{equation}\label{intrep2}
\phi(s)=2zs-\ln\,s-\ln(1+s).
\end{equation}
When $z>0$ the contour  $\cal L$ is a vertical line with $\Re s>0$; when $z$ is complex we can deform the contour in order to keep convergence. For the same purpose  we can rotate the path of integration for $P_\nu^\mu(z)$ in \eqref{intrep1}.

For certain choices of $p_\mu(s)$ and $q_\mu(s)$ these integrals give representations of the functions $Y_n^{\mu}(z)$, $F_n^{\mu}(z)$, and $U_n^{\mu}(z)$. We have
\begin{equation}\label{intrep3}
Y_n^{\mu}(\zeta)=\frac{(2\nu z)^{n+\mu+1}}{\Gamma(n+\mu+1)} P_\nu^\mu(z), \quad p_\mu(s)=\frac{s^{\mu}}{\sqrt{s(1+s)}},
\end{equation}
\begin{equation}\label{intrep4}
U_n^{\mu}(\zeta)=\frac{(-1)^n(2\nu z)^{n+\mu+1}e^{-2\nu z}}{\Gamma(n+\mu+1)}P_\nu^\mu(z), \quad p_\mu(s)=\frac{(1+s)^\mu}{\sqrt{s(1+s)}},
\end{equation}
\begin{equation}\label{intrep5}
F_n^{\mu}(\zeta)=\frac{n!}{(2\nu z)^{n}} Q_\nu^\mu(z), \quad q_\mu(s)=\frac{(1+s)^{-\mu}}{\sqrt{s(1+s)}}.
\end{equation}
The multi-valued functions in $\phi(s)$, $p_\mu(s)$, and $q_\mu(s)$ have there principal branches and are real for $s>0$.

The representations in \eqref{intrep3} and \eqref{intrep4} follow from the well-known integral
\begin{equation}\label{intrep6}
U(a,c,z)=\frac{1}{\Gamma(a)}\int_0^\infty t^{a-1}(1+t)^{c-a-1} e^{-zt}\,dt, \quad \Re\,a,z>0,
\end{equation}
the first line in \eqref{ykum}, and the second line in \eqref{ykumnegFU}. 
For \eqref{intrep5} we refer to the first line in \eqref{ykumnegFU} and the integral representation (see \cite[p.~46]{Slater:1960:CHF})
\begin{equation}\label{intrep7}
\frac{1}{\Gamma(c)}{}_1F_1(a;c;z)=\frac{z^{1-c-}e^z}{2\pi i}\int_\calL e^{zs}(1+s)^{a-c}s^{-a} \,ds, 
\end{equation}
where $\calL$ is a vertical line with $\Re s>0$.

\subsection{Construction of the expansions}\label{sec:deriv}

We use the saddle point method for obtaining asymptotic expansions of the integrals in \eqref{intrep1}. The saddle points follow from
the equation $\phi^\prime(s)=0$, where
\begin{equation}\label{der1}
\phi^\prime(s)=\frac{2z s^2+2(z-1)s-1}{s(1+s)},
\end{equation}
and are given by
 \begin{equation}\label{der2}
{s_{+}=\frac{1-z+\sqrt{1+z^2}}{2z},} \quad
{s_{-}=\frac{1-z-\sqrt{1+z^2}}{2z}.}
\end{equation}
When $z>0$ the saddle points are well-separated, with $-1<s_-< -\frac12$ and $s_+>0$. We have the following limits: $\lim_{z\to0} s_{+}=+\infty$ and  $\lim_{z\to\infty}s_{+}=0$.

Also,
\begin{equation}\label{der3}
s_{+}(1+s_{+})=\frac{1+\sqrt{1+z^2}}{2z^2},
\end{equation}
and
\begin{equation}\label{der4}
\phi(s_+)=1-z+\ln(2z)+\eta, \quad 
\phi^{\prime\prime}(s_+)=\frac{4z^2\sqrt{1+z^2}}{1+\sqrt{1+z^2}},
\end{equation}
with $\eta$ defined in \eqref{etat0}.

We use Laplace's method with the transformation 
\begin{equation}\label{der5}
\phi(s)-\phi(s_+)=\tfrac12\phi^{\prime\prime}(s_+) w^2, \quad \sign(w)=\sign(s-s_+).
\end{equation}
We have $s=w+\bigO(w^2)$ as $w\to0$. The integrals in \eqref{intrep1} become
\renewcommand{\arraystretch}{1.75}
\begin{equation}\label{der6}
\begin{array}{l}
\dsp{P_\nu^\mu(z)=e^{-\nu\phi(s_+)}\int_{-\infty}^\infty f(w)\, e^{-\tfrac12\nu\phi^{\prime\prime}(s_+) w^2}\,dw,\quad f(w)=p_\mu(s)\frac{ds}{dw},}\\
\dsp{Q_\nu^\mu(z)=\frac{e^{\nu\phi(s_+)}}{2\pi i}\int_{-i\infty}^{i\infty}  g(w)\, e^{\tfrac12\nu\phi^{\prime\prime}(s_+) w^2}\,dw,\quad g(w)=q_\mu(s)\frac{ds}{dw}.}
\end{array}
\renewcommand{\arraystretch}{1}
\end{equation}

By expanding $f(w)=\sum_{k=0}^\infty f_k w^k$ and  $g(w)=\sum_{k=0}^\infty g_k w^k$ we obtain the asymptotic expansions
 \renewcommand{\arraystretch}{1.75}
\begin{equation}\label{der7}
\begin{array}{l}
\dsp{P_\nu^\mu(z)\sim f_0 e^{-\nu\phi(s_+)}\sqrt{\frac{2\pi }{\nu \phi^{\prime\prime}(s_+)}}
\sum_{k=0}^\infty \frac{F_{k}(\mu,z)} {\nu^k},}\\
\dsp{Q_\nu^\mu(z)\sim g_0 \frac{e^{\nu\phi(s_+)}}{2\pi}\sqrt{\frac{2\pi }{\nu \phi^{\prime\prime}(s_+)}}
\sum_{k=0}^\infty (-1)^k \frac{G_{k}(\mu,z)} {\nu^k},}
\end{array}
\renewcommand{\arraystretch}{1}
\end{equation}
where (see also \eqref{der3})
 \renewcommand{\arraystretch}{1.75}
\begin{equation}\label{der8}
\begin{array}{l}
\dsp{F_k(\mu,z)= \frac{(\frac12)_k 2^k }{(\phi^{\prime\prime}(s_+))^k}\frac{f_{2k}}{f_0}}, \quad
\dsp{f_0=p_\mu(s_+)},\\
\dsp{G_k(\mu,z)= \frac{(\frac12)_k 2^k }{(\phi^{\prime\prime}(s_+))^k}\frac{g_{2k}}{g_0}}, \quad
\dsp{g_0=q_\mu(s_+)},
\end{array}
\renewcommand{\arraystretch}{1}
\end{equation}
because $ds/dw=1$ at $w=0$. 

By using \eqref{der2}--\eqref{der4} and \eqref{intrep4}--\eqref{intrep5} it follows that
 \begin{equation}\label{der9}
Y_n^{\mu}(\zeta)\sim\frac{\left(1-z+\sqrt{1+z^2}\right)^{\mu}\sqrt{z}}{(1+z^2)^{\frac14}}
\frac{e^{\nu z-\nu \eta}}{\Gamma^*(\nu+\mu+\frac12)}
\sum_{k=0}^\infty \frac{F_{k}^{(1)}(\mu,z)} {\nu^k},
\end{equation}
 \begin{equation}\label{der10}
U_n^{\mu}(\zeta)\sim(-1)^n\frac{\left(1+z+\sqrt{1+z^2}\right)^{\mu}\sqrt{z}}{(1+z^2)^{\frac14}}
\frac{e^{-\nu z-\nu \eta}}{\Gamma^*(\nu+\mu+\frac12)}
\sum_{k=0}^\infty \frac{F_{k}^{(2)}(\mu,z)} {\nu^k},
\end{equation}
\begin{equation}\label{der11}
F_n^{\mu}(\zeta)\sim\Gamma^*\left(\nu+\tfrac12\right)\frac{\left(1+z-\sqrt{1+z^2}\right)^\mu\sqrt{z}}{(1+z^2)^{1/4}}e^{-\nu z +\nu\eta}\sum_{k=0}^\infty(-1)^k\frac{G_k(\mu,z)}{\nu^k}.
\end{equation}
The coefficients $F_{k}^{(1)}(\mu,z)$ are obtained from \eqref{der8} and the function $f(w)$ of \eqref{der6} with the function 
$p_\mu(s)$ as given in  \eqref{intrep3}, and $F_{k}^{(2)}(\mu,z)$ follow from taking the function $p_\mu(s)$ as given in  \eqref{intrep4}. The function $\Gamma^*$ is the slowly varying part of the corresponding gamma function. That is,
\begin{equation}\label{der12}
\Gamma^*(\nu+\alpha)=\frac{\Gamma(\nu+\alpha)}{\sqrt{2\pi}\,\nu^{\nu+\alpha-\frac12}e^{-\nu}}\sim
\sum_{k=0}^{\infty}\frac{\gamma_k(\alpha)}{\nu^k}, \quad \gamma_0(\alpha)=1, \quad  \nu\to\infty.
\end{equation}
The coefficients  $\gamma_k(\alpha)$ follow from standard methods for the gamma function;  see 
\S\ref{subsec:coeffY}.

The final form of the expansion of $Y_n^{\mu}(\zeta)$ given in \eqref{Ynp15} in Theorem~\ref{thm2} can be obtained by dividing the expansion in \eqref{der9} by the expansion of $\Gamma^*(\nu+\mu+\frac12)$ given in \eqref{der12}, and similar for the other expansions. This gives for $k=0,1,2,\ldots$
\begin{equation}\label{der13}
A_{k}(\mu,z)=F_{k}^{(1)}(\mu,z)-\sum_{j=0}^{k-1} A_{j}(\mu,z)\gamma_{k-j}\left(\mu+\tfrac12\right),
 \end{equation}
\begin{equation}\label{der14}
B_{k}(\mu,z)=F_{k}^{(2)}(\mu,z)-\sum_{j=0}^{k-1} B_{j}(\mu,z)\gamma_{k-j}\left(\mu+\tfrac12\right), 
\end{equation}
\begin{equation}\label{der15}
C_{k}(\mu,z)=\sum_{j=0}^k (-1)^jG_{j}(\mu,z)\gamma_{k-j}\left(\tfrac12\right).
 \end{equation}

\subsection{Computation of the coefficients}\label{subsec:coeffY}
To compute the coefficients $F_k(\mu,z)$ and $A_k(\mu,z)$ we need the coefficients in the expansion
\begin{equation}\label{comp1}
s=s_++\sum_{k=1}^\infty s_kw^k,
\end{equation}
which follow from \eqref{der5}. We write, as in  \eqref{etat}, 
 \begin{equation}\label{comp2}
 t=\frac{1}{\sqrt{1+z^2}}
 \end{equation}
 and obtain
\renewcommand{\arraystretch}{2}
 \begin{equation}\label{comp3}
\begin{array}{l}
\dsp{s_1=1,\quad s_2=\frac{2-t}{6}, \quad s_3=\frac{(1-t)(5t^3-6t^2+2)}{18t^2},} \\
\dsp{s_4=-\frac{z(1-t)(40t^4-65t^3+24t^2-2t+4)}{135t^2}.}
\end{array}
\end{equation}
\renewcommand{\arraystretch}{1.0}%
With these coefficients we can compute the coefficients $f(w)$ and $g(w)$ of \eqref{der6} by choosing the appropriate $p_\mu(s)$ and $q_\mu(s)$. 

To obtain the coefficients in \eqref{der13}--\eqref{der15} we first compute $\gamma_k(\mu+\frac12)$ that appear in \eqref{der12}. We have
\renewcommand{\arraystretch}{1.75}
\begin{equation}\label{comp4}
\begin{array}{l}
\dsp{\gamma_0(\mu+\tfrac12)=1, \quad \gamma_1(\mu+\tfrac12)=\tfrac{1}{24}\left(-1+12\mu^2\right),}\\
\dsp{\gamma_2=\tfrac{1}{1152}\left(1+48\mu-24\mu^2-192\mu^3+144\mu^4\right),}\\
\dsp{\gamma_3=\tfrac{1}{414720}\left(1003-720\mu-17100\mu^2+11520\mu^3+32400\mu^4-}
\right.\\
\quad\quad\quad\quad
\left.
\dsp{34560\mu^5+8640\mu^6\right),}\\
\dsp{\gamma_4=\tfrac{1}{39813120}\left(-4027-288864\mu+151824\mu^2+1618560\mu^3-}
\right.\\
\quad
\left.
\dsp{1239840\mu^4-1645056\mu^5+2177280\mu^6-829440\mu^7+103680\mu^8\right).}
\end{array}
\end{equation}
\renewcommand{\arraystretch}{1}%

\begin{figure}
\caption{\small 
Saddle point contours of the first  integral in  \eqref{intrep1} for $z= e^{i\theta}$, $\theta=k\pi/10$, $k=1,2,3,4,5$.
\label{BP.fig1}}
\begin{center}
\epsfxsize=10cm \epsfbox{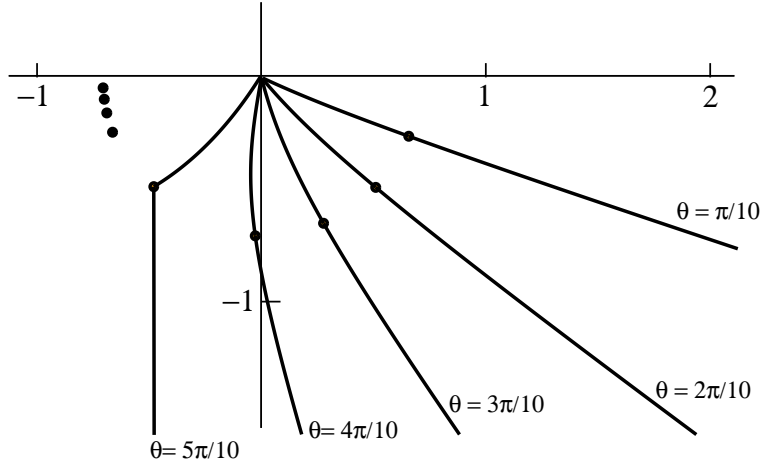}
\end{center}
\end{figure}

\subsection{Extending the result to complex values of  \protectbold{z}}\label{subsec:extz}

From \cite[p.~378]{Olver:1997:ASF}  it follows that the expansions in \eqref{Knuzuniform} and \eqref{Ynp17} hold for large values of $\nu$ and are uniformly valid for complex values of $z$ inside the sector $\vert\phase\,z\vert\le \frac12\pi-\delta$ with $\delta$ a small positive number. 
As can be seen from the front factor and the coefficients, it becomes invalid when $z$ approaches~$\pm i$. In that case the singularities of the functions $f$ and $g$ in \eqref{der6} approach the origin. 

The singularities come from those of the mapping in \eqref{der5}. This mapping does not depend on $\mu$ and, hence, Laplace's method remains applicable for all fixed values of $\mu$, and also \eqref{Ynp15} is uniformly valid for complex values of $z$ inside the sector $\vert\phase\,z\vert\le \frac12\pi-\delta$.
 
For complex $z$ inside the sector $\vert\phase\,z\vert< \frac12\pi$ the saddle points given in \eqref{der2} move into the complex plane, and it is for all these values of $z$ possible to find a single saddle point contour from $0$ through $s_+$ such that $\phase(zs)=0$ at infinity. If $\vert\phase\,z\vert\le \frac12\pi-\delta$ the singular points of the mapping in \eqref{der5} and of the function $f$ and $g$ in \eqref{der6}  stay away from the origin.

In Figure~\ref{BP.fig1} we show the saddle points $s_+$ (black balls) and the corresponding saddle point contours of the first integral in \eqref{intrep1} for $z= e^{i\theta}$, $\theta=k\pi/10$, $k=1,2,3,4,5$. The black balls at the left are the saddle points $s_-$. When $z=i$ the saddle points $s_+$ and $s_-$ coincide at $s=-\frac12(1+i)$. The saddle point contours of the second integral in \eqref{intrep1} are the paths of steepest ascent  of the first integral.

\subsection{Comparison with earlier expansions}\label{subsec:comp}
In \cite{Dunster:2001:UAE} the expansions of $Y_n^\mu(z)$ are given for large $n$ with possibly large values of $\mu$ as well. This makes comparison with our expansions rather complicated. In fact in Dunster's expansions given in \cite[\S\S6-7]{Dunster:2001:UAE}  the expansions can be re-expanded for small values of a parameter $\alpha$ corresponding with (in our notation) $\mu/\nu$, and in this way our results of the present section may be obtained. Dunster has used Olver's theory \cite{Olver:1997:ASF} for linear differential equations of second order, with bounds for the remainders in the expansions and a recursion formula for the coefficients.

Starting from an integral we show how to include $\mu$  as a second large parameter, leaving out the details. 
Let $\mu =\alpha\nu$ and write \eqref{Ynp2} in the form 
\begin{equation}\label{D1}
Y_n^{\mu}(\zeta)=\frac{(2\nu z)^{n+\mu+1}}{\Gamma(n+\mu+1)} 
\int_0^\infty \frac{e^{-\nu\psi(s)}}{\sqrt{s(1+s)}}\,ds,
\end{equation}
where $\zeta=1/(\nu z)$, $\nu=n+\frac12$, and
\begin{equation}\label{D2}
\psi(s)=2zs-(1+\alpha)\ln\,s-\ln(1+s).
\end{equation}
Then the saddle point analysis for the result of \S\ref{sec:Beszpos} can be repeated, giving an expansion that holds again in the sector $\vert\phase\,z\vert\le \frac12\pi-\delta$  and $\alpha\ge-1+\varepsilon$, with $\delta, \varepsilon$  small positive numbers. In a similar way the expansions of \S\ref{sec:Beszneg} can be modified.

With respect to the results in terms of elementary functions given in \cite{Wong:1997:AEG} we observe the following points.
\begin{itemize}
\item
A remarkable point is that  for $Y_n^\mu(z)$ the point $z=0$  excluded, whereas in our results this point is accepted as long we approach it inside the sectors $\vert\phase(\pm z)\vert \le \frac12\pi-\delta$. In fact the results are essentially the same as our results, although the notation and scaling of the parameters is different. For example, the factor $e^{-1/z}$ at the left-hand side of \cite[(2.26)]{Wong:1997:AEG} can be combined with $e^{(n+1)f(\zeta_+,\alpha)}$ to give a regular expression as $z\to0$. Also, the domains of validity are different and are simpler in our case (just the sectors $\vert\phase(\pm z)\vert \le \frac12\pi-\delta$). 
\item
The expansions are derived from the integral representation \eqref{Besint}, after a transformation. A detailed discussion is given about the location of saddle points and paths of steepest descent and the domains for the expansions to be derived. Our start is from two integrals related to those for the Kummer functions. 
\item
The reader has to become familiar with several aspects of the detailed description of the paths of integration, the domains of validity and the proper choices of the branches of some multi-valued quantities. In this sense, our approach is more accessible for the reader who wants to use the results and to construct more terms.
\item
Complex quantities arise in the coefficients and front factors (which, of course for real $z$ will provide real expansions). Our expansions  show quantities that are real for real $z$.
\item
It is not indicated how the results reduce to the well-known expansions of the modified Bessel functions; relations with the Kummer functions are not given.
\end{itemize}

\section{Expansions in terms of modified Bessel functions}\label{sec:gen}
The expansions for $Y_n^\mu(1/(\nu z))$ in the previous section \S\ref{sec:elem} become invalid when $z$ approaches the points $\pm i$, because in that case the saddle points coincide. As shown in \cite{Dunster:2001:UAE} and \cite{Wong:1997:AEG} it is possible to derive uniform expansions in terms of Airy functions, and these expansions are valid in large $z-$domains.

For the modified Bessel functions $I_\nu(\nu z)$ and  $K_\nu(\nu z)$ similar asymptotic phenomena arise when $z$ approaches the points $\pm i$, and the expansion in terms of Airy functions is available in the literature. In fact expansions for the Hankel functions and the ordinary Bessel functions can be used. See \cite[]{Olver:2010:BFS} and \cite[Chapter~11]{Olver:1997:ASF}. 

Because the asymptotic phenomena of the generalized Bessel polynomials $Y_n^\mu(z)$ for large $n$ and fixed $\mu$ are the same as those of the polynomial $Y_n^0(z)$, we approach the problem for obtaining uniform expansions by expanding the generalized polynomials in terms of the modified Bessel functions $K_\nu(z)$ (with $\nu=n+\frac12$), which are the same as the reduced Bessel polynomials $Y_n^0(z)$ (see \eqref{ybes}). By using the existing results for the Bessel functions a complete description is available in this way.

We summarize the results of this section as follows.

\begin{theorem}\label{thm4}
For $n\to \infty$ we have the asymptotic expansion
\begin{equation}\label{gen13}
Y_n^{\mu}(\zeta)\sim\frac{(2\nu z)^{\mu} n!\,e^{\nu z}}{\Gamma(n+\mu+1)}\sqrt{\frac{2\nu z}{\pi}}\,
\left(K_\nu(\nu z)\sum_{k=0}^{\infty}\frac{C_k}{\nu^k}+
K_\nu^\prime(\nu z)\sum_{k=0}^{\infty}\frac{D_k}{\nu^k}\right),
\end{equation}
and the expansion holds uniformly with respect to all $z$. Here, $\zeta=1/(\nu z)$,  $\nu=n+\frac12$,
\renewcommand{\arraystretch}{1.75}
\begin{equation}\label{gen17}
C_0=\left(2^{-\frac12}e^{-\frac34\pi i}\right)^{\mu},\quad
C_1=\tfrac1{24}(1-i)\mu(\mu-1)(-2\mu+1+3i)C_0,
\end{equation}
and
\begin{equation}\label{gen18}
D_0=\tfrac12(1-i)\mu C_0,\quad
D_1=-\tfrac1{24}i \mu^2 (\mu-1) (-\mu+2+3 i)  C_0.
\end{equation}
\renewcommand{\arraystretch}{1.0}
\end{theorem}

\subsection{The construction of the expansion}\label{subsec:constr}
To start the construction of the expansion we  write \eqref{yint} in the form
\begin{equation}\label{Ynp2}
Y_n^{\mu}(\zeta)=\frac{(2\nu z)^{n+\mu+1}}{\Gamma(n+\mu+1)} 
\int_0^\infty \frac{s^{\mu}}{\sqrt{s(1+s)}}e^{-\nu\phi(s)}\,ds,
\end{equation}
where again
\begin{equation}\label{Ynp3}
\zeta=\frac{1}{\nu z},\quad  \nu=n+\tfrac12, \quad \phi(s)=2zs-\ln\,s-\ln(1+s).
\end{equation}
We write
\begin{equation}\label{gen01}
f_0(s)= s^\mu= A_0+B_0s +\phi^{\prime}(s)g_0(s),
\end{equation}
and substitute $s=s_+$ and $s=s_-$ to obtain
\begin{equation}\label{gen02}
A_0=\frac{s_+f_0(s_-)-s_-f_0(s_+)}{s_+-s_-},\quad
B_0=\frac{f_0(s_+)-f_0(s_-)}{s_+-s_-}.
\end{equation}
Putting \eqref{gen01} into \eqref{Ynp2} we obtain
\begin{equation}\label{gen03}
Y_n^{\mu}(\zeta)=A_0\Phi_0+B_0\Phi_1+ \frac{(2\nu z)^{n+\mu+1}}{\Gamma(n+\mu+1)}
\int_0^\infty \frac{\phi^{\prime}(s)g_0(s)}{\sqrt{s(1+s)}}e^{-\nu\phi(s)}\,ds,
\end{equation}
where
\begin{equation}\label{gen04}
\Phi_0=\frac{(2\nu z)^{\mu} n!}{\Gamma(n+\mu+1)}Y_n^0(\zeta),\quad
\Phi_1=\frac{(2\nu z)^{\mu-1} (n+1)!}{\Gamma(n+\mu+1)}Y_n^1(\zeta).
\end{equation}
By using \eqref{ybes} and \eqref{ydiff} it follows that
\renewcommand{\arraystretch}{2.0}
\begin{equation}\label{gen05}
\begin{array}{l}
\dsp{\Phi_0=\frac{(2\nu z)^{\mu} n!}{\Gamma(n+\mu+1)}\sqrt{\frac{2\nu z}{\pi}}\,
e^{\nu z}K_{\nu}(\nu z),} \\
\dsp{\Phi_1=\frac{(2\nu z)^{\mu} n!}{2\Gamma(n+\mu+1)}
\sqrt{\frac{2\nu z}{\pi}}\,e^{\nu z}
\left((1/z-1)K_\nu(\nu z)-K_\nu^\prime(\nu z)\right)}.
\end{array}
\end{equation}
\renewcommand{\arraystretch}{1.0}%
In the second line we can also write \cite[p.~234]{Temme:1996:SFI}
\begin{equation}\label{gen06}
(1/z-1)K_\nu(\nu z)-K_\nu^\prime(\nu z)=K_{\nu+1}(\nu z)-K_\nu(\nu z),
\end{equation}
but we prefer the notation with the derivative, because the asymptotic expansions of $K_\nu(\nu z)$ and $K_\nu^\prime(\nu z)$ are quite related and usually presented together. 

The next step is to use integration by parts in \eqref{gen03}, and this gives
\begin{equation}\label{gen07}
Y_n^{\mu}(\zeta)=A_0\Phi_0+B_0\Phi_1+ \frac{(2\nu z)^{n+\mu+1}}{\nu\Gamma(n+\mu+1)}
\int_0^\infty\frac{ f_1(s)}{\sqrt{s(1+s)}}e^{-\nu\phi(s)}\,ds,
\end{equation}
where
\begin{equation}\label{gen08}
f_1(s)=\sqrt{s(1+s)}\frac{d}{ds}\frac{g_0(s)}{\sqrt{s(1+s)}}.
\end{equation}
Repeating this procedure by writing for $k\ge0$
\begin{equation}\label{gen09}
f_k(s)= A_k+B_ks +\phi^{\prime}(s)g_k(s),\quad f_0(s)=s^\mu,
\end{equation}
\begin{equation}\label{gen10}
A_k=\frac{s_+f_k(s_-)-s_-f_k(s_+)}{s_+-s_-},\quad
B_k=\frac{f_k(s_+)-f_k(s_-)}{s_+-s_-},
\end{equation}
\begin{equation}\label{gen11}
f_{k+1}(s)=\sqrt{s(1+s)}\frac{d}{ds}\frac{g_k(s)}{\sqrt{s(1+s)}}
=g_k^\prime(s)-\frac{2s+1}{2s(s+1)}g_k(s),
\end{equation}
we obtain for $K\ge0$
\renewcommand{\arraystretch}{2}
\begin{equation}\label{gen12}
\begin{array}{l}
\dsp{Y_n^{\mu}(\zeta)=\Phi_0\sum_{k=0}^{K-1}\frac{A_k}{\nu^k}+
\Phi_1\sum_{k=0}^{K-1}\frac{B_k}{\nu^k}\,+}\\
\quad\quad\quad\quad\quad\quad
\dsp{ \frac{(2\nu z)^{n+\mu+1}}{\nu^M\Gamma(n+\mu+1)}
\int_0^\infty\frac{ f_K(s)}{\sqrt{s(1+s)}}e^{-\nu\phi(s)}\,ds.}
\end{array}
\end{equation}
\renewcommand{\arraystretch}{1.0}

We rearrange the expansion by using \eqref{gen05} and writing
\begin{equation}\label{gen14}
C_k= A_k+\frac{1-z}{2z}B_k, \quad D_k=-\tfrac12B_k, \quad k=0,1,2,\ldots\,
\end{equation}
to obtain the expansion given in Theorem~\ref{thm4}.

\begin{remark}\label{Genrem1}
{\rm
To compute the coefficients $A_k, B_k$ defined in \eqref{gen10} and the functions $f_k(s)$, say, by using a computer algebra package, it is convenient to write the functions $f_k(s)$ in the form of two-point Taylor expansions at the saddle points $s_+$ and $s_-$. More details on this method can be found in \cite{Lopez:2002:TPT}, \cite{Lopez:2004:MPT}, \cite{Vidunas:2002:SEC}.
}
\end{remark}

\begin{remark}\label{Genrem2}
{\rm
For integer values of $\mu$ we have the following simple cases.
\begin{enumerate}
\item
For $\mu=0, 1,2,\ldots$ the expansion in \eqref{gen12} has a finite number of terms which can also be obtained from the recursion in \eqref{yrecmu}. 
\item
For $\mu=-1,-2,-3,\ldots$ we can also obtain an exact result.  When $\mu=-1$ we have
\begin{equation}\label{gen19}
C_k=\frac{z-1}{2^k}, \quad D_k=-\frac{z}{2^k},\quad k=0,1,2,\ldots,
\end{equation}
and we can sum the convergent series when $2\nu =2n+1>1$. This gives a result that corresponds to the relation in \eqref{yrecmu} with $\mu=-1$.
\end{enumerate}
}
\end{remark}

\section{Concluding remarks}\label{Conc}
In \S\ref{sec:Besalt} we have given a new simple expansion of $Y_n^{\mu}(z)$  that is valid outside a compact neighborhood of the origin in the $z-$plane and new forms of expansions in terms of elementary functions valid  in the sectors $\vert\phase(\pm z)\vert \le \frac12\pi-\delta$ not containing the turning points $z=\pm i/n$. To avoid mappings for obtaining expansions in terms of Airy functions we have given expansions in terms of modified Bessel functions. For these functions very detailed Airy-type expansions are available, which can be used to obtain similar expansions of $Y_n^{\mu}(z)$. 

\section{Appendix: Expansions of the modified Bessel functions}\label{sec:modBes}

Because we compare the expansions for the Bessel polynomials to those of the modified Bessel functions, we summarize a few details about the uniform expansions of the $K-$ and $I-$Bessel functions.

We have \cite{Olver:2010:BFS}\footnote{http://dlmf.nist.gov/10.41},  \cite[p.~378]{Olver:1997:ASF} 
\begin{equation}\label{Knuzuniform}
K_{\nu}(\nu z)\sim \sqrt{\frac{\pi}{2\nu}}\frac{e^{-\nu\eta}}{(1+z^2)^{1/4}}
\sum_{k=0}^{\infty}(-1)^k\frac{u_k(t)}{\nu^k},
\end{equation}
\begin{equation}\label{Inuzuniform}
I_{\nu}(\nu z)\sim\frac{1}{ \sqrt{2\pi\nu}}\frac{e^{\nu\eta}}{(1+z^2)^{1/4}}
\sum_{k=0}^{\infty}\frac{u_k(t)}{\nu^k},
\end{equation}
which hold when $\nu\to\infty$, uniformly with respect to $z$ such that 
$\vert\phase\, z\vert\le\frac12\pi-\delta$, $\delta$ being an arbitrary positive number in $(0,\frac12\pi)$. Here,
\begin{equation}\label{etat}
t=\frac{1}{\sqrt{1+z^2}}, \qquad
\eta=\sqrt{1+z^2}+\log\frac{z}{1+\sqrt{1+z^2}}.
\end{equation}

The first coefficients $u_k(t)$ are
\begin{equation}\label{ukt}
u_0(t)=1, \qquad u_1(t)=\frac{3t-5t^3}{24}, 
\qquad u_2(t)=\frac{81t^2-462t^4+385t^6}{1152},
\end{equation}
and other coefficients can be obtained by applying the formula
\begin{equation}\label{recuk}
u_{k+1}(t)=\tfrac 12 t^2(1-t^2)u_k'(t)+\tfrac 18 \int_0^t (1-5s^2)u_k(s) ds,
\quad k=0,1,2,\ldots.
\end{equation}

For the derivatives we have
\begin{equation}\label{Knuzduniform}
K_{\nu}^\prime(\nu z)\sim - \sqrt{\frac{\pi}{2\nu}}\frac{(1+z^2)^{1/4}}{z}e^{-\nu\eta}
\sum_{k=0}^{\infty}(-1)^k\frac{v_k(t)}{\nu^k},
\end{equation}
\begin{equation}\label{Inuzduniform}
I_{\nu}^\prime(\nu z)\sim\frac{1}{ \sqrt{2\pi\nu}}\frac{(1+z^2)^{1/4}}{z}e^{\nu\eta}
\sum_{k=0}^{\infty}\frac{u_k(t)}{\nu^k},
\end{equation}
where
\begin{equation}\label{vkt}
v_0(t)=1, \qquad v_1(t)=\frac{-9t+7t^3}{24}, 
\qquad v_2(t)=\frac{-135t^2+594t^4-455t^6}{1152},
\end{equation}
and other coefficients can be obtained by applying the formula
\begin{equation}\label{recvk}
v_{k}(t)=u_{k}(t)+t(t^2-1)\left(\tfrac12u_{k-1}(t)+tu_{k-1}^\prime(t)\right),
\quad k=0,1,2,\ldots.
\end{equation}

The expansions in \eqref{Knuzuniform} and \eqref{Inuzuniform} become invalid when $z$ approaches the turning points $\pm i$.  In that case expansions are available  in terms of Airy functions. First the functions $I_{\nu}(z)$ and $K_{\nu}(z)$ should be written in terms of ordinary Bessel functions, and then the results for these functions can be used; see \cite[p.~419--426]{Olver:1997:ASF}.

\section*{Acknowledgments}
The authors thank the referee for helpful comments on the first version of the paper. 
The authors acknowledge financial support from {\it Gobierno of Navarra}, Res.~07/05/2008; 
NMT acknowledges support from {\emph{Ministerio de Ciencia e Innovaci\'on}},
project MTM2009--11686.



\end{document}